\documentclass[12pt]{article}

\usepackage{amsmath}
\usepackage{amssymb}
\usepackage{amsthm}

\usepackage[dvips]{graphicx}
\usepackage{color}
\usepackage{epsfig}
\usepackage[mathscr]{eucal}

\newtheorem{theorem}{Theorem}

\newtheorem{prop}[theorem]{Proposition}

\newtheorem{cor}[theorem]{Corollary}

\newtheorem{principle}{Principle}

\newtheorem{conjecture}[principle]{Conjecture}

\newtheorem{problem}[principle]{Problem}

\newcommand{\df}{\ensuremath{\textrm{def}}}
\newcommand{\floor}[1]{\ensuremath{\left \lfloor {#1} \right \rfloor}}
\newcommand{\ceil}[1]{\ensuremath{\left \lceil {#1} \right \rceil}}

\newcommand{\trip}[1]{\ensuremath{\textrm{Trip} \! \left( {#1} \right )}}

\newcommand{\Z}{{\mathbb Z}}

\newcommand{\cE}{{\cal E}}

\newcommand{\cG}{{\cal G}}

\allowdisplaybreaks

\title{Collinear Points in Permutations}
\author{Joshua N. Cooper \\ \small Courant Institute of Mathematics \\ \small New York University, New York, NY \\ \\ J\'{o}zsef Solymosi \\ \small Department of Mathematics \\ \small University of British Columbia, Vancouver, BC}

\begin{document}

\maketitle

\begin{abstract} Consider the following problem: how many collinear triples of points must a transversal of $\Z_n \times \Z_n$ have?  This question is connected with venerable issues in discrete geometry.  We show that the answer, for $n$ prime, is between $(n-1)/4$ and $(n-1)/2$, and consider an analogous question for collinear quadruples.  We conjecture that the upper bound is the truth and suggest several other interesting problems in this area.
\end{abstract}

In \cite{R51}, Erd\H{o}s offered a construction concerning the ``Heilbronn Problem''.  What is the smallest $A$ so that, for any choice of $n$ points in the unit square, some triangle formed by three of the points has area at most $A$?  His elegant construction of a point-set with large minimum-area triangle ($\sim n^{-2}$) is as follows: take the smallest prime $p \geq n$, and let the set of points be $\{p^{-1}(x,x^2 \pmod{p}): x \in \Z_p\}$.  (If necessary, throw out a few points so that there are $n$ left.) It is easy to see that this set has no three collinear points, and therefore any three points form a nondegenerate lattice triangle -- which must have area at least $p^{-2}/2 \gg n^{-2}$.

Another area in which collinear triples of points on a lattice arise is in connection with the so-called ``no-three-in-line'' problem, dating back at least to 1917 (\cite{D17}).  Is it possible to choose $2n$ points on the $n$-by-$n$ grid so that no three are collinear?  Clearly, if this is the case, then $2n$ is best possible.  Guy and Kelly (\cite{GK68}) conjecture that, for sufficiently large $n$, not only is it true that every set of $2n$ points has a collinear triple, but that it is possible to avoid collinear triples in a set of size $(\alpha - \epsilon)n$ and impossible to avoid them in a set of size $(\alpha + \epsilon)n$, where $\alpha = (2\pi^2/3)^{1/3} \approx 1.874$ and $\epsilon > 0$.

In this note, we address the question of when it is possible to avoid collinear triples {\it modulo} $n$, particularly in the case of transversals (i.e., graphs of permutations) and when $n$ is prime.

\section{Results}

Suppose that $S$ is a subset of $\Z_n \times \Z_n$.  We say that a set of points $X \subset S$ is ``collinear'' if there are parameters $a,b,c \in \Z_n$ so that each $(i,j) \in X$ lies on the line $\{(x,y) \in \Z_n \times \Z_n : ax + by = c\}$.  If $n$ is prime, the ``slope'' of a such a line is defined to be $ab^{-1}$ if $b \neq 0$ and $\infty$ otherwise.  Clearly, each pair of points is collinear, and we say that that slope of a pair of points is the slope of the line containing them.  (It is easy to see that this is well-defined.)

If $f : \Z_n \rightarrow \Z_n$ is some function, then we say that the ``graph'' of $f$ is the set $\{(x,f(x)):x \in \Z_n\}$.  Often we will say that a set of points of a {\it function} are collinear, where really we mean that a set of points of its graph are.

We are interested in the number of collinear triples in subsets $S$ of $\Z_n \times \Z_n$, and conditions guaranteeing that there is at least one.  Note that, if $n$ is prime, then if $|S| = m$, then $(m-1)(n-1)+1 \leq n^2$.  Indeed, fix some point $s \in S$; then each line through $s$ contains $n-1$ points of $\Z_n \times \Z_n$, plus $s$ itself.  If we count each line passing through $s$ and some other $t \in S$, then the total number of points is $(m-1)(n-1)+1$, since the lines are pairwise disjoint except for $s$.  Therefore the inequality holds, and moreover, $m \leq n+2$.  If some point $s = (s_1,s_2)$ has the property that $(t,s_2) \not \in S$ for any $t \neq s_1$, then $(m-1)(n-1)+n \leq n^2$, and therefore $m \leq n+1$.  Finally, if $s$ has the property that $(t_1,s_2) \not \in S$ and $(s_1,t_2) \not \in S$ for any $t_1,t_2 \neq s_1$, then this count yields $(m-1)(n-1)+n+(n-1) \leq n^2$, i.e., $m \leq n$.  Therefore, in general, (1) a set with no collinear triples can have at most $n+2$ points, (2) a set with no collinear triples which has some column (or row) containing only one point can have at most $n+1$ points, and (3) a set with no collinear triples which has a point lying in an otherwise empty column {\it and} row can have at most $n$ points.  It is interesting to ask, then, how many collinear triples must a permutation of $\Z_n$ have?  Define $\Psi(n)$ to be the minimum number of collinear triples in any permutation of $\Z_n$.

\subsection{Permutations}

\begin{theorem} \label{onetriplethm} Suppose $\sigma$ is a permutation of $\Z_n$ for $n>2$ prime.  Then $\sigma$ contains at least one collinear triple of points, i.e., $\Psi(n) > 0$.
\end{theorem}
\begin{proof} Suppose not.  Note that for every pair of points in the graph of $\sigma$, the slope of that pair must be in $1, \ldots, n-1$.  Partition the $\binom{n}{2}$ pairs into classes according to their slopes.  Since there are $n-1$ classes, at least $n/2$ pairs lie in some class.  Each of these pairs lies on some line of the same slope, and no two of them lie on the same line, since then we would a collinear triple.  Therefore, we have at least $n/2$ lines, each of which contains two points.  However, $n$ is odd, so there are in fact at least $(n+1)/2$ disjoint lines containing two points.  This amounts to $n+1$ points in total, a contradiction.
\end{proof}

By the argument above, it is clear that the $\binom{n}{2}$ pairs lie in at most $(n-1)^2/2 = \binom{n}{2} - (n-1)/2$ lines.  Therefore, a natural question is to ask, if $K$ pairs are assigned to $L$ families, what is the minimum number of triples occurring entirely within some family?  That is, suppose $\cG$ is a graph with $K$ edges, and $E(\cG)$ is partitioned into $L$ ``lines'' $\cE_1,\ldots,\cE_L$.  What is the least possible value of
$$
\sum_{i=1}^L \binom{|\bigcup \cE_i|}{3},
$$
over all $\cG$ and partitions $\{\cE_i\}$?  Call this quantity $T(K,L)$.  If some $\cG$ and a partition of its edges achieves this bound, then the edges which belong to a given line must span the minimum number of vertices.  That is, if $|\cE_i| = m_i$, then
$$
\left|\bigcup \cE_i \right| = \ceil{\frac{1+\sqrt{1 + 8m_i}}{2}} =_\df \tau(m_i),
$$
because this is the smallest number $k$ so that $\binom{k}{2} \geq m_i$.  Therefore, $T(K,L)$ is the least possible value of
$$
\trip{m_1,\ldots,m_L} =_\df \sum_{i=1}^L \binom{\tau(m_i)}{3}
$$
over all $L$-tuples $m_1,\ldots,m_L$ of nonnegative integers whose sum is $K$.\\

Write $\rho(t)$ for $t - 2 \floor{t/2}$, i.e., the parity of $t$.

\begin{prop} \label{KLprop} For $K \leq 3L$, $T(K,L) = \max\{\ceil{(K-L)/2},0\}$.
\end{prop}
\begin{proof} If $K \leq L$, it is clear that we may set $m_i = 1$ for $1 \leq i \leq K$ and $m_i = 0$ otherwise, so that $T(K,L)=0$.  Therefore, suppose $L < K \leq 3L$.  We make the following claim: there is some partition $K = m_1 + \cdots + m_L$ minimizing $\trip{m_1,\ldots,m_L}$ which has $\floor{(K-L)/2}$ indices $i$ so that $m_i = 3$, $\rho(K-L)$ indices $i$ so that $m_i = 2$, $K - 3 \floor{(K-L)/2} - 2 \rho(K-L)$ indices $i$ so that $m_i = 1$, and the other $m_i = 0$.  We prove the claim by induction on $K$.  The base case $K=L$ we dealt with above.  If $K = L+1$ or $K = L+2$ then we may set $m_1 = K - L + 1$ and $m_i = 1$ for $2 \leq i \leq L$, resulting in $\trip{m_1,\ldots,m_L} = 3$, which is clearly best possible.  Suppose, then, that $K \geq L+3$.

By the pigeonhole principle, we may assume that $m_1 > 1$ without loss of generality.  Furthermore, since $\tau(2)=\tau(3)$ and increasing $m_1$ can only decrease
$$
\trip{m_2,\ldots,m_L} = \trip{m_1,\ldots,m_L} - \binom{\tau(m_1)}{3},
$$
we have that $m_1 \geq 3$ unless $K \leq 2$, i.e., we are in the base case.  If $m_1 = 3$, then
$$
\trip{m_1,\ldots,m_L} = 1 + \trip{m_2,\ldots,m_L} = 1 + \ceil{\frac{K-3-(L-1)}{2}} = \ceil{\frac{K-L}{2}}.
$$
If, on the other hand, $m_1 > 3$, then $\tau(m_1) \geq 4$ and $L \geq 2$.  In that case,
\begin{align*}
\trip{m_1,\ldots,m_L} & = \binom{\tau(m_1)}{3} + \trip{m_2,\ldots,m_L} \\
& = \binom{\tau(m_1)}{3} + \ceil{\frac{K-L-m_1+1}{2}} \\
& \geq \frac{K-L+1}{2} + \binom{m_1}{3}  - \frac{m_1}{2} \\
& > \ceil{\frac{K-L}{2}},
\end{align*}
where the second line follows from induction and the fourth from the fact that $\binom{\tau(m_1)}{3}-\frac{m_1}{2} > 0$ for $m_1 > 3$.  Since this contradicts the minimality of the partition, we must have $m_1 = 3$.  Furthermore, there are
$$
1 + \floor{\frac{K-L-2}{2}} = \floor{\frac{K-L}{2}}
$$
indices $i$ with $m_i$ = 3, $\rho(K - L - 2) = \rho(K-L)$ indices $i$ with $m_i$ = 2,
$$
K - 3 - 3 \floor{(K-L-2)/2} - 2 \rho(K-L-2) = K - 3 \floor{(K-L)/2} - 2 \rho(K-L)
$$
indices $i$ with $m_i = 1$ and the rest zeroes.

\end{proof}

We may apply this immediately to the original question by setting $K = \binom{n}{2}$ and $L = (n-1)^2/2$.

\begin{cor} For $n>2$ prime, $\Psi(n) \geq \ceil{(n-1)/4}$.
\end{cor}

\noindent On the other hand, we have the following.

\begin{prop} \label{upperbound} For $n>2$ prime, $\Psi(n) \leq (n-1)/2$.
\end{prop}
\begin{proof} Define the function $f:\Z_n \rightarrow \Z_n$ by $f(x)=x^{-1}$ if $x \neq 0$ and $f(0)=0$.  This is clearly a permutation, and we show that its graph has exactly $(n-1)/2$ collinear triples.  First, suppose that $1 \leq x < y < z \leq n-1$.  It is easy to see that the condition that $(x,f(x))$, $(y,f(y))$, and $(z,f(z))$ are collinear is equivalent to the statement that
\begin{align*}
&(f(z)-f(x))(y-x)-(f(y)-f(x))(z-x) = \\
&(z^{-1}-x^{-1})(y-x)-(y^{-1}-x^{-1})(z-x) = 0 \pmod{n}.
\end{align*}
Multiplying by $xyz$ yields
$$
y(x-z)(y-x)-z(x-y)(z-x) = (y-z)(x-z)(y-x) = 0,
$$
which is impossible because $x$, $y$, and $z$ are distinct.  Therefore, if $f$ exhibits a collinear triple, it must have some point with abscissa $0$.  Without loss of generality, we may assume that $x=0$.  Therefore, if $0=x<y<z\leq n-1$,
$$
z^{-1}y-y^{-1}z = 0,
$$
so that $y^2 = z^2$, i.e., $y = \pm z$.  Since $y \neq z$, $y=-z$, and there are $(p-1)/2$ unordered triples of the type $\{0,y,-y\}$.
\end{proof}

In fact, we believe the following to be true.

\begin{conjecture} $\Psi(n) = (n-1)/2$ for $n>2$ prime.
\end{conjecture}

We have included computational data supporting this conjecture in Table 1.

\begin{table}[!ht]
\begin{center}
\begin{tabular}{|c|c||c|c||c|c|} \hline
$x$ & $\Psi(x)$ & $x$ & $\Psi(x)$ & $x$ & $\Psi(x)$ \\
\hline \hline 1  & 0 & 7 & 3 & 13 & 6 \\
\hline 2         & 0 & 8 & 0 & 14 & 9 \\
\hline 3         & 1 & 9 & 5 & 15 & 6 \\
\hline 4         & 0 & 10 & 2 & 16 & 4 \\
\hline 5         & 2 & 11 & 5 & 17 & 8 \\
\hline 6         & 0 & 12 & 0 & 18 & $\leq 16$ \\
\hline
\end{tabular}
\caption{The first $18$ values of $\Psi(x)$.}
\end{center}
\end{table}

\begin{problem} What about $n$ composite?
\end{problem}

About this question, unfortunately, we can say nothing.\\

Note that, by the proof of Proposition \ref{upperbound}, any fractional linear transformation of the form $(ax+b)/(cx+d)$ with $c \neq 0$ (along with $-dc^{-1} \mapsto ac^{-1}$) gives rise to $(n-1)/2$ collinear triples for $n > 2$ prime.

\begin{conjecture} The function
$$
g(x) = \left \{ \begin{array}{ll} x/(x-1) & \textrm{if $x \neq 1$} \\ 1 & \textrm{if $x = 1$} \end{array} \right .
$$
is the lexicographic-least permutation with $(n-1)/2$ collinear triples for $n>2$ prime.
\end{conjecture}

\subsection{Quadruples}

The permutation from Proposition \ref{upperbound} has the property that, if we remove the point $(0,0)$, the resulting graph is a collinear triple-free set of $n-1$ points with no two on a single row or column -- showing that Theorem \ref{onetriplethm} is tight.  A moment's reflection also reveals that it has no collinear 4-tuple.

There are permutations with many collinear triples which have no collinear 4-tuples, however.  Consider $h(x) = x^3$, a function on $\Z_n$, where $n > 2$ is prime and congruent to $2$ mod $3$.  Then $h$ is a permutation, since the unique solution of $x^3=c \pmod{n}$ is $x=c^{(n-1)/3} \pmod{n}$.  Clearly, $x^3-ax-b=0$ has at most three solutions, so no line intersects the graph of $h$ (or $h_0$) in four points.  Furthermore, $x^3-ax-b$ cannot have exactly two roots in $\Z_n$ unless it has a double root, since
$$
x^3 - ax - b = (x-r_1)(x-r_2)(x-r_3)
$$
implies that $r_3 = br_1^{-1}r_2^{-1}$ unless $r_1 = 0$ or $r_2 = 0$.  If $r_1 = 0$, then $x^3-ax-b = x^3 - x^2 (r_2+r_3) + x r_2r_3$, so $r_3 = -r_2$, and either both $r_2$ and $r_3$ are in $\Z_n$ or neither is.  (The same holds for $r_1$ and $r_3$ if $r_2 = 0$.)  If it has a double root, then
$$
x^3 - ax - b = (x-r_1)^2(x-r_2) = x^3 - (2r_1+r_2)x^2+(r_1^2+2r_1r_2)x-r_1^2r_2,
$$
so $r_2 = -2r_1$.  Therefore, $a = 3r_1^2$ and $b = -2r_1^3$.  This gives $n-1$ lines containing exactly two points ($r_1 = 0$ gives a triple root).

Each of the $X$ collinear triples of $h$ contains three pairs of points.  No pair is contained in two triples, since two collinear triples that intersect in two points form a collinear quadruple.  Furthermore, only $n-1$ pairs are not contained in {\it some} collinear triple, so we have $3 X + (n-1) = \binom{n}{2}$, or $X = (n-1)(n-2)/6$.  We may conclude that, for each prime $n > 2$ congruent to $2$ mod $3$, there is some permutation with $(n-1)(n-2)/6$ collinear triples, but no collinear quadruples.  Unfortunately, such a construction cannot work for $p = 1 \pmod{3}$, since no cubic permutation polynomials exist for such $p$ (q.v. \cite{LN97}).

Now, suppose that a permutation has $X$ collinear triples but no collinear quadruples.  Each of the $X$ triples contains three pairs of points and, again, each such pair appearing in a triple appears in only one of them.  Therefore, we must have $3 X \leq n(n-1)/2$, i.e., $X \leq \floor{n(n-1)/6}$ -- and this holds for {\it any} $n$.

By the above observations, $\limsup_{n \rightarrow \infty} CT(n)/n^2 = 1/6$.  However, the question remains what the lower bound is when $n \not \equiv 2 \pmod{3}$, or when $n$ is composite.  There is also a gap of about $n/3$ between the upper and lower bounds even in the case of $n \equiv 2 \pmod{3}$.  We also wish to know:

\begin{problem} What is the maximum number $CT_0(n)$ of collinear triples in a subset of $\Z_n \times \Z_n$ which has no collinear quadruples?
\end{problem}

\subsection{Pair Packing}

We return to the argument of Proposition \ref{KLprop}.  Consider the following ``greedy'' process for finding an assignment of $K$ pairs to $L$ lines.  We proceed through the lines one at a time, and when we reach $\cE_i$ with $K^\prime$ pairs/edges unaccounted for, we place into it
$$
\min\left \{K^\prime, \binom{\tau(\ceil{K^\prime/(L-i+1)})}{2} \right \}
$$
edges, arranged as a graph on $\tau(\ceil{K^\prime/(L-i+1)})$ vertices.  That is, we distribute the remaining edges as equally as possible into the remaining lines, but round the number of edges we place into the current line $\cE_i$ up to the nearest triangular number whenever possible.  One might conjecture that this process results in the optimal configuration, i.e, minimizing $T(K,L)$ -- but it does not.  Indeed, already for $K = 28$ and $L=2$, the optimal configurations are $(m_1,m_2) = (21,5)$ and $(20,6)$, neither of which has a line with $m_i = 15 = \binom{\tau(28/2)}{2}$.

\begin{problem} \label{prob3} Describe those configurations which achieve $T(K,L)$.
\end{problem}

It is easy to see that, for each $K$ and $L$, {\it some} optimal configuration has at most one $m_i$ which is not a triangular number.  Furthermore, if $m_i \leq r$ for some $r$, then we may move $q$ points from the largest line $m_j = \binom{a}{2} + q$, $1 \leq q \leq a$, to line $i$ and in the process change the number of total triples by
$$
\binom{\tau(q+r)}{3} - \binom{a+1}{3} + \binom{a}{3} \leq \binom{\tau(a+r)}{3} - \binom{a}{2}.
$$
Therefore, we have a contradiction if $\binom{\tau(a+r)}{3} < \binom{a}{2}$.  Since $\tau(x) \leq 2 + \sqrt{2x}$,
$$
\binom{\tau(a+r)}{3} \leq \frac{(2 + \sqrt{2(a+r)})^3}{6}
$$
For $r$ large, $a = 2 r^{3/4}$ provides a contradiction.  Therefore,

\begin{prop} If $\min_i m_i = r$, $r$ sufficiently large, in a configuration $(m_1,\ldots,m_L)$ achieving $T(K,L)$, then $\max_i m_i \leq 2 r^{3/2}$.
\end{prop}

It is easy to see that the function
$$
g(x) = \binom{(1+\sqrt{1+8x})/2}{3} = \frac{x}{6} \left (\sqrt{1+8x}-3 \right)
$$
is concave-up for $x \geq 0$.  Therefore, by Jensen's Inequality, $T(K,L) \geq L g(K/L)$.  This, in particular, implies that the number of collinear triples in a subset of $\Z_n \times \Z_n$ of cardinality $X$, $n$ prime, is $\gg X^3/n^2$ when $X/n \rightarrow \infty$.

Surely, more than this can be said concerning Problem \ref{prob3}.

\section{Acknowledgements}

Thank you to Greg Martin, Jozef Skokan and Joel Spencer for
valuable discussions and insights.

\end{document}